\def\note{April 1996.
 This is an update of my problem list \cite{m31}.
 I tried to include as many references as I could think of.
 If you know anything about these problems or could supply
 any missing references or corrections (missing
 attributions or misattributions), please let me know.
 }

\documentstyle[12pt]{article}
\input amssym.def
\input amssym

\def\reals{{\Bbb R}}
\def\csp{2^{\omega}}
\def\bsp{\omega^{\omega}}
\def\infsets{[\omega]^{\omega}}
\def\rect{{\cal R}}

\def\answer#1{\par\medskip \begin{quote} Answer: #1 \end{quote}\medskip}
\def\where{\marginpar{where???}}

\newcount\probno
\def\topic#1{\subprobno=1\advance\probno by 1
 \newpage\begin{center} {\the\probno . #1}  \end{center}}
\def\topicx#1{\subprobno=1\advance\probno by 1
 \par\bigskip\begin{center} {\the\probno . #1}  \end{center}}

\newcount\subprobno
\def\prob{\par\medskip\noindent
 {\kern -.5in \parbox{.45in}{\the\probno .\the\subprobno}}
 \advance\subprobno by 1}

\begin{document}

\begin{center}
                    Some interesting problems\\
                   Arnold W. Miller\footnote{\note}
\end{center}

\topicx{Analytic sets}

\prob (Mauldin) Is there a $\Sigma^1_1$ set X universal for
$\Sigma^1_1$ sets which are not Borel?  Suppose $B\in\Sigma^1_1$ and
for every  Borel A, $A\leq_W B$.  Does this imply that for every
$\Sigma^1_1$ A, $A\leq_W B$. (This refers to Wadge reducible.)
\answer{The first question was answered by Hjorth \cite{hjorth1}
who showed that it is independent.}

\prob A subset $A\subset\bsp$ is compactly-$\Gamma$ iff for every
compact $K\subset\bsp$ we have that $A\cap K$ is in $\Gamma$.
Is it consistent relative to ZFC that compactly-$\Sigma^1_1$ implies
$\Sigma^1_1$?  (see Miller-Kunen \cite{m17}, Becker \cite{beckanal})

\prob (\cite{m17})
Does $\Delta^1_1$ = compactly-$\Delta^1_1$ imply
$\Sigma^1_1$ = compactly-$\Sigma^1_1$?

\prob (Prikry see \cite{fried}) Can $L\cap\bsp$ be a nontrivial $\Sigma^1_1$
set?  Can there be a nontrivial perfect set of constructible reals?
\answer{This is almost settled by Velickovic-Woodin \cite{vw}.}

\prob (Ostaszewski \cite{ostem})
Consider Telgarsky's game $G(T)$ where $T\subseteq \csp$.
Player I plays a countable cover of $T$
Player II chooses one- say  $X_n$.
Player I wins iff $\cap\{cl(X_n):n\in\omega\}\subseteq T$.
It is known that
\par (a) Player I has winning strategy iff T is analytic.
\par (b) If there exists $A$ an analytic subset of $cl(T)$ not Borel separated
from $T$, then Player II has a winning strategy.
\par\noindent Is the converse of (b) true?

\prob Does there exists an analytic set which is not
Borel modulo Ramsey null?  Same question for the ideal generated
by closed measure zero sets.

\prob (Sierpinski \cite{sier70}) Does there exists
an analytic set of reals $E$ such that every (uncountable) analytic
set of reals is the one-to-one continuous image of $E$?

\prob (Jockusch 10-95) Let for $A\subseteq \omega$ let
$D(A)=\{a-b:a,b\in A\}$.  Is the set
$\{D(A): A\subseteq\omega\}$ Borel?

\prob Suppose $I$ is a $\sigma$-ideal generated by
its $\Pi^0_2$ members.
Then is it true that for any analytic set $A$
either $A\in I$  or $A$ contains a $\Pi^0_3$ set not in I?
This is suggested by a theorem of Solecki \cite{solec}
that says that for any
$\sigma$-ideal $I$ generated its closed members and analytic set $A$,
either $A\in I$  or $A$ contains a $G_\delta$ set not in I.

\topic{Axiom of Determinacy}

\prob Does AD imply that $2^{\omega_1}$ is the $\omega_1$ union of meager
sets? \answer{Yes, Becker \cite{becker}.}

\prob Does AD imply that there does not exist $\omega_2$ distinct
$\Sigma_2^1$ sets?

\prob Is there a hierarchy of $\Delta^1_2$ sets?

\prob Does AD imply every set is Ramsey?
\answer{Yes, if also assume $V=L[\reals]$ for references
see Kanamori \cite{kan} page 382.}

\prob (V. Delfino \cite{cab}) (Conjecture) If $f:\csp\mapsto\csp$ is
Turing invariant ($x\equiv_T y\to f(x)\equiv_T f(y)$) then
there exists $z$ such that either for every
$x\geq_{T} z\;\;f(x)\geq_{T} x$ or there exist $c$ such that
for every $x\geq_{T} z\;\;f(x)\equiv_T c$.

\topic{Combinatorial cardinals less than the continuum}

\prob (van Douwen \cite{vandou})
If every $\omega_2$ descending sequence in
$P(\omega)/{\mbox{finite}}$
has something beneath it  is it true that every family of $\omega_2$
sets with the IFIP has something beneath it? (does $\goth t=\goth p$).

\prob (Hechler \cite{hech}) Let M be a countable transitive model of ZFC.
Does there exists $M[f_n : n\in\omega]$ a generic extension with each
$f_n\in\bsp$ such that for all $n\in\omega$
$f_n$ eventually dominates every element of
$$M[f_m : m>n]\cap\bsp.$$
For something similar with Sacks forcing see Groszek \cite{gros}.

\prob (Dow) Does the following imply $p(\kappa)$:
$\forall X,Y\subset\infsets,|X|,|Y|<\kappa$ and $\forall A\in X,B\in Y$
$A\cap B$ finite, there exists $U\subset\omega$ such that
for all $V\in (X\cup Y)$
\centerline{ $(U\cap V)$ is infinite iff $V\in X$  }

\prob  Can the least $\kappa$ such that Indep($\kappa$) fails
have cofinality $\omega$?  Indep($\kappa$) means that  every
family B of $\kappa$ infinite subsets of $\omega$ there exists
an infinite subset Z of $\omega$ such that for every $A\in B$,
$|Z\cap A|=|Z\setminus A|=\omega$.

\prob (Kunen \cite{kunma}) Let $\goth m$ be the smallest cardinal for which
MA$_{\goth m}$ fails.  Can we have $\omega_2=cof(\goth m)<\goth m$?

\prob (Scheepers 7-91) Is it consistent that $\aleph_\omega$
embeds into $(\bsp, \leq^*)$
but not $\aleph_{\omega+1}$?

\topic{MAD families}

\prob (Roitman) Is it consistent that
every maximal almost disjoint
family in $\infsets$ has cardinality greater than  $\omega_1$, but
there exists a dominating family F in $\bsp$ of cardinality
$\omega_1$? For a related result see Shelah \cite{shelboulder}.

\prob (van Douwen) CH implies there exist $F\subset\bsp$ which
is maximal with respect to eventually different functions
which is also maximal
with respect to infinite partial functions also.
Is there always such a
one?  What is the cardinality of the smallest?

\prob (Cook, Watson) Consider paths in $\omega\times\omega$. CH implies
there is a MAD family of paths.  Is there always one?
\answer{No, Steprans \cite{step},
still open for dimensions $\geq 3$}

\prob (Milliken) A maximal almost disjoint family X is a separating
family iff for all $Q\in\infsets$
$$\{Q\cap P{\;\;|\;\;} P\in X\}$$
has size continuum or is finite.
Are there always separating families?

\prob (Erd\"{o}s, Hechler \cite{erdhech})
Does MA plus the continuum is larger than $\aleph_{\omega+2}$ imply
that there is no mad family on $\aleph_{\omega}$ of size
$\aleph_{\omega+1}$?

\prob (Kunen) Call $I\subset\infsets$ an independent splitting family
if I is independent ( every finite boolean combination of elements of I
is infinite) and splitting ( for every $f:I\mapsto 2$ there does not
exist an infinite X such that for every $A\in I$, $X\subset^* A^{f(A)}$,
where $A^0=A$ and $A^1=\omega\setminus A$.)
If CH or MA then there does exists an independent splitting family.
In ZFC is there one?

\prob (Fleissner) If there is a Luzin set, then is there a MAD family
of size $\omega_1$?

\topic{Forcing}

\prob (S. Friedman, R. David) Let $P_n=2^{<\omega_n}$.  Does
forcing with $\Pi_{n\in\omega}P_n$ add a Cohen
 subset of $\omega_{\omega+1}$?
\answer{Yes, Shelah \cite{shel1}}

\prob (Kunen) Force with perfect $P\subset\csp$ such that
for every $I\in\infsets$ $\pi_I:P\mapsto 2^I$ does not have a countable
range.  Is $\omega_1$ collapsed?

\prob (van Douwen, Fleissner) Is it consistent with not
CH that for P a c.c.c partial order of
size continuum there exists a sequence $G_{\alpha}$ for $\alpha<\omega_1$
of P-filters such that for every dense set $D\subseteq P$ all but countably
many of the $G_{\alpha}$ meet D.
\answer{No, Todorcevic \cite{tod}}

\prob Is there a Truss-like characterization of eventually different reals?
How about infinitely equal reals?  (see Truss \cite{truss} for
characterization of dominating reals).

\prob (Kunen \cite{kun}) Does there exists an $\omega_1$ saturated
$\sigma$-ideal
in the Borel subset of $\csp$ which is invariant under homeomorphisms
induced by permutations of $\omega$ and different from the meager ideal,
measure zero ideal, and the intersection ideal?
\answer{Partial answer Kechris-Solecki \cite{kecsol}}

\prob (van Mill) Is it consistent that every c.c.c. boolean algebra
which can be embedded
into $P(\omega)/{\mbox{finite}}$ is $\sigma$-centered?

\prob Suppose $M\subseteq M[f]$ are models of ZFC and
for every $g\in\bsp\cap M$ there exists infinitely
many  $n\in\omega$ such that  $g(n)=f(n)$.
Must there exists a real $x\in M[f]$ which is Cohen over
$M$?  (If there are two such infinitely equal reals, then
there must be a Cohen real, see \cite{m11} and \cite{slalom}.)

\topic{Measure theory}

\prob (Mauldin, Grzegorek) Is it consistent that the continuum is RVM
and all sets of reals of cardinality $\omega_2$ have zero measure?
\answer{No, apparently from the Gitik-Shelah
Theorem (see Fremlin \cite{fremreal} 6F) it was deduced by
Prikry and Solovay that if $\kappa$ is
real-valued measurable, then there are Sierpinski sets of all cardinalities
less than $\kappa$.}

\prob (Fremlin) Can the cardinality of the
least cover of the real line by measure zero
sets have countable cofinality?
\answer{Yes, Shelah \cite{shel2}}

\prob (Erd\"{o}s) For every sequence converging to zero does there exist
a set of positive measure which does not contain a similar sequence?
Falconer \cite{falc}
has shown that if the sequence converges slowly enough there
does exist such a set of positive measure.  H.I. Miller \cite{him}
has shown the analogous statement for Baire category to be false.
I showed that for every sequence there exist a partition
of the reals into two sets neither of which contains a sequence similar
to the given one.

\prob (Erd\"{o}s) Suppose for every $n\in\omega$ the set
$A\cap [n,\infty)$ has positive measure.  Must $A$ contain arbitrarily
long arithmetic progressions?

\answer{Several mathematicians have pointed out this is trivial.
Probably I misquoted Erd\"{o}s.  I scribbled it down after one
of his talks when the universe was younger.
To quote Just \cite{just13-1}:
``The answer to 6.4 seems to be trivially
`yes', unless you want the differences to be integers; then the answer seems
to be trivially `no', unless you want the measure to be positive in
EVERY interval, in which case
the answer may not be so trivial. So, what should the
problem really look like?''
}

\prob Is it possible to have a Loeb-Sierpinski set of cardinality greater
than $\omega_1$? See Leth-Keisler-Kunen Miller \cite{m25} and
Miller \cite{m27}.

\prob (Louveau) If a subset A of the plane has positive measure
and contains the diagonal, then
does there exist a set B in the line of positive outer measure such than
$B^2$ is a subset of A?
\answer{According to Burke \cite{bur}, Fremlin and Shelah proved
this fails in the Cohen real model.}

\topic{Borel hierarchies}

\prob Is it consistent that for every countable ordinal $\alpha$ there
exists a $\Pi^1_1$ set of Baire order $\alpha$?
See Miller \cite{m1}.

\prob Is it consistent that for every uncountable separable metric space
X there exists a X-projective set not Borel in X? See Miller
\cite{m9},\cite{m28}.

\prob Is it consistent that the set of all Baire orders is the same
as the set of even ordinals $\leq\omega_1$? See Miller \cite{m33}.

\prob Is it true that if X is a $Q_{\alpha}$-set and Y is a  $Q_{\beta}$-set
and $2\leq\alpha<\beta$ then $|X|<|Y|$? \cite{m1}

\prob Does $\rect^{\omega_2}_{\omega_1}=P(\omega_2\times\omega_2)$ and
$2^{\omega}=\omega_2$ imply that
$2^{\omega_1}=\omega_2$?   ($\rect^{\omega_2}$ is the family of abstract
rectangles in $\omega_2\times\omega_2$ and the lower subscript is the
level of the Borel hierarchy.)

\prob Does $\rect^{\omega_2}_{\omega}=P(\omega_2\times\omega_2)$ imply that
for some $n<\omega\; \rect^{\omega_2}_n=P(\omega_2\times\omega_2)$?

\prob Does $|X|=\omega_1$ imply that X is not a $Q_{\omega}$-set?

\prob (Mauldin) Is it consistent that there exists a separable metric space
X of Baire order less than $\omega_1$ (i.e. for some $\alpha<\omega_1$
every Borel subset of X is $\Sigma^0_{\alpha}$ in $X$) but not every
relatively analytic set is relatively Borel?

\prob Can the Borel hierarchy on cubes in $\rect^3$ behave differently than
the Borel hierarchy on rectangles in $\rect^2$?

\prob (Ulam \cite{ulam}) Is there a separable metric space of
each projective
class order? ($\Sigma^1_2$-forcing?) See Miller \cite{m9},\cite{m28}.

\prob In the Cohen real model is there an uncountable separable metric
space of Baire order 2?  In the random real model are there any separable
metric spaces of Baire order between 2 and $\omega_1$?
\answer{Answered in Miller \cite{m33}.}

\prob What can we say about hierarchy orders involving
difference hierarchies or even abstract $\omega$-boolean operations?

\prob (Stone) Is it consistent to have a  Borel map $f:X\mapsto Y$ where
$X$ and  $Y$ are metric spaces and $f$ has the property that there is
no bound less than $\omega_1$ on the Borel complexity of $f^{-1}(U)$ for
$U\subseteq Y$ open?  Fleissner \cite{fleiss}
shows that it is consistent there is no such $f$ using
a supercompact.

\prob (Ciesielski-Galvin \cite{cg}) Let $P_2(\kappa)$ be the family of
all cylinder
sets in $\kappa^3$ (where cylinder means $A\times B$ where $A\subseteq \kappa$
and $B\subseteq \kappa^2$ or anything that could be obtained like this by
permuting the three coordinates.) Is it consistent that the $\sigma$-algebra
generated by $P_2(\goth c^{++})$ is equal to all subsets of
$(\goth c^{++})^3$?

\prob (Ciesielski)  Suppose every subset of $\omega_2\times\omega_2$
is in the $\sigma$-algebra generated by the abstract rectangles.  Does this
continue to hold after adding $\omega_1$-Cohen reals?

\prob (Fleissner \cite{fleisq})
If $X$ is a Q-set of size $\omega_1$ then is $X^2$ a Q-set?
(Not necessarily true if for $X$ of cardinality $\omega_2$.)

\topic{Involving $\omega_1$}

\prob (Jech-Prikry \cite{jp}) Is it consistent that
there exists a family $F\subset\omega^{\omega_1}$
of cardinality less than $2^{\omega_1}$, such that for every
$g\in\omega^{\omega_1}$ there exist $f\in F$ such that for every
$\alpha<\omega_1, g(\alpha)<f(\alpha)$.

\prob (Frankiewicz \cite{fran}) Is it consistent that
 $\beta\omega\setminus\omega$ is homeomorphic
to $\beta\omega_1\setminus\omega_1$?

\prob Is it consistent to have CH, $2^{\omega_1}>\omega_2$, and
there exists $F\subset[\omega_1]^{\omega_1}$ of cardinality
$\omega_2$ such that for every $A\subset\omega_1$ there exists
$B\in F$ such that $B\subseteq A$ or $B\cap A=\emptyset$?

\prob (Kunen) Is it consistent to have  $2^{\omega_1}>\omega_2$ and
there exists $F\subset[\omega_1]^{\omega_1}$ of cardinality
$\omega_2$ such that for every uncountable $A\subset\omega_1$ there exists
$B\in F$ such that $B\subseteq A$?

\prob (Kunen) Is it consistent to have  a uniform ultrafilter
on ${\omega_1}$ which is generated by fewer than $2^{\omega_1}$ sets?

\prob (Prikry \cite{m20}) Is it consistent there exists an
$\omega_1$ generated ideal J such that $P(\omega_1)=P(\omega)/J$?

\prob (Comer) If C and D are homeomorphic to $2^{\omega_1}$ then
is $C\cup D$? (Say if both are subsets of $2^{\omega_1}$.)

\prob (Nyikos) If $C\times D$ is homeomorphic to $2^{\omega_1}$ then
must either C or D be homeomorphic to $2^{\omega_1}$?

\prob (CH) Let n(X) be the cardinality of the smallest family of meager sets
which cover X.  Can the cof(n($(2^{\omega_1})_{\delta}$))
($G_{\delta}$-topology) be $\omega$ or $\omega_1$?

\prob (Velickovic) Is it consistent that every Aronszajn line L
contains a Countryman type?

\prob Does PFA imply that any two Aronszajn types contain uncountable
isomorphic subtypes?

\topic{Set theoretic topology}

\prob Is it consistent to have no P-points or Q-points?  A P-point
is an ultrafilter U on $\omega$ with the property that every
function $f:\omega\mapsto\omega$ is either constant or finite-to-one
on an element of U.
A Q-point
is an ultrafilter U on $\omega$ with the property that every finite-to-one
function $f:\omega\mapsto\omega$ is one-to-one
on an element of U.
Shelah \cite{shel3} showed it is consistent there are no P-points
and Miller \cite{m2} showed that
it is consistent there are no Q-points.  Roitman and Taylor
showed that if the continuum is $\leq \omega_2$, then there must
be a P-point or a Q-point.

\prob (M.E. Rudin) Is there always a small Dowker space?
\answer{Yes?, Balogh \cite{baldowk}, Kojman-Shelah \cite{kojshe}}

\prob (Charlie Mills) In infinite dimensional Hilbert space is a sphere
coverable by fewer than continuum other spheres?

\prob Is it consistent that $\omega^{\omega_1}$ is pseudonormal?
(Pseudonormal means disjoint closed sets can be separated if
at least one is countable.)

\prob (van Douwen) Is it consistent to have
$c(U(\omega_1))<d(U(\omega_1)$?  (c is celluarity, d is density and
U is uniform ultrafilters.)

\prob  Is the box product of countably many copies of the unit interval
coverable by countably many zero dimensional sets?

\prob (Hansell) Is there a non-zero-dimensional Q-space?
Can there be
a nonzero dimensional metric space in which every subset is $G_{\delta}$?

\prob (Bing) Suppose $D_n$ a subset of the plane is homeomorphic to a disk
and for every $n\in\omega$ $D_{n+1}\subseteq D_{n}$, then does
$\cap_{n\in\omega}D_n$ have the fixed point property?

\prob (Sikorski see \cite{kur})
Does there exist two closed 0-dimensional nonhomeomorphic subsets
of the plane such that each is homeomorphic to an open subset of the other?

\prob (Ancel 11-93)
Is there a separable Hausdorff space in which every basis has cardinality
$2^{2^c}$? (c = the cardinality of the set of real numbers.)

\prob (Gulko 1995) Is there a model of ZFC in which there is a maximal
almost family $M$ on $\omega$ such that for any
point $x\in \beta\omega\setminus\omega$ there exists a countable
subfamily of $M$ such that $x$ is in its closure.


\topic{Model Theory}

\prob (Vaught) Does every countable first order theory have countably many
or continuum many countable models up to isomorphism?  How about for universal
theories of a partial order?  See Becker \cite{beck}.

\prob (Martin) Show that if T is a countable first order theory with
fewer than continuum many countable models up to isomorphism, then
every countable model of T has an isomorphism class which is at most
$\Sigma^0_{\omega+\omega+1}$.

\prob (Caicedo \cite{caic}) Does every theory in $L_{\omega_1,\omega}$ have an
independent axiomatization?

\prob Does V=L imply there exists a complete theory T such that
$$\{\alpha : L_{\alpha}\models T\}$$ is an unbounded subset of $\omega_1$?
\answer{Yes, Hjorth \cite{hjorem}.}

\prob (Miller \cite{m13})
Are there any properly $\Sigma^0_{\lambda+1}$ isomorphism classes
for $\lambda$ a countable limit ordinal?
\answer{Yes, Hjorth \cite{hjorthclass}}

\prob Is there a theory with exactly $\omega_1$ rigid countable models up to
isomorphism? (same for minimal models)
\answer{Yes, for minimal models, Hjorth \cite{hjorthmin}.
Yes, for rigid models, Hjorth \cite{hjorthrig}}

\prob (Baldwin) Are there continuum many complete $\omega$-stable theories
in a finite language?  (I don't know if this is still open?)

\prob (Miller-Manevitz \cite{m18})
Is it consistent that there exists a model of ZFC, M, such the
unit interval of M is Lindelof and $\omega^M$ is $\omega_1$-like?

\prob (Miller \cite{m13})
Does there exists a pseudo-elementary class in the language of
one unary operation with exactly $\omega_1$ nonisomorphic countable
models?  (There is pseudo $L_{\omega_1,\omega}$ class.)

\prob (Mati Rubin) Does there exists an embedding of the rationals into
themselves such that no between function is elementarily extendable?

\prob Duplicate of 10.6

\prob (suggested by Fuhrken \cite{fuh} \cite{ef})
Can we have a model with exactly one undefinable $L_{\omega_1,\omega}$
element?

\prob
Can we have a complete first order theory $T$ with models of size
$\aleph_{2n}$ for $n< \omega$ (but not of size $\aleph_{2n+1}$) and
$\aleph_\omega < c$?

\topic{Special subsets of the real line}

\prob (Mauldin, Grzegorek) Is it consistent that every
universally measurable set has the Baire property?  See
Corazza \cite{cor}.

\prob (Mauldin) Are there always $>\goth c$ universally measurable sets?
(same for restricted Baire property)

\prob (Galvin) Does every Sierpinski set have strong first category?
\answer{Bartoszynski-Judah \cite{bj} showed that it is consistently yes.
Yes, Pawlikowski \cite{paw}.}

\prob (Galvin, Carlson) Is the union of two strong first category
sets a set of strong first category?
\answer{Not necessarily, Bartoszynski-Shelah \cite{barshel}}

\prob Does there exists a perfectly meager $X\subseteq\reals^n$ which is not
zero-dimensional?

\prob (Kunen) Is it consistent that for every uncountable
$X\subseteq \reals$
there exists a measure zero set M such that $X+M$ has positive
outer measure? See Erdos-Kunen-Mauldin \cite{ekm}.

\prob (Sierpi\'{n}ski \cite{sier}) A set of reals X is a J-set iff for
uncountable
$Y\subseteq X$ there exist a perfect $P\subseteq X$ such that $P\cap Y$
is uncountable.  If we assume CH, then a set is a J-set iff it is
$\sigma$-compact.  Is it consistent with not CH that
J-set $ = \sigma$-compact?

\prob (Fremlin-Miller \cite{m22})
Is there always an uncountable subset of the reals which is hereditary
with respect to property M?
\answer{Yes, see \cite{m34}}

\prob  Consider the three  non c.c.c ideals: $(s)_0$-sets, Ramsey null sets,
and $\sigma$-compact sets.
What can one say about the properties of add, cov, non,
and cof?  ( add = additivity of ideal, cov = smallest cardinality of a
cover of reals by subset of ideal, non = smallest cardinality of a
set of reals not in ideal,
cof = cofinality = smallest cardinality of a family of sets in ideal
which has the property
that every set in ideal is covered by some element of the family.)

\prob Consider the notion of Laver null sets.  This is defined analogously
to Ramsey null sets, but use Laver forcing instead of Mathias forcing.
The analogue of Galvin-Prikry Theorem is true here.  What other results
also go thru?  What ideals arise from other notions of forcing?
What about Silver forcing?
What notions of forcing arise from infinite combinatorial theorems?
(For example, Carlson's infinite version of the Hales-Jewett theorem
\cite{carl}.)

\prob (Judah-Shelah \cite{js}) Is the Borel conjecture plus the existence of
a Q-set consistent?

\prob (Daniel, Gruenhage).  Given a set of reals $X$ and ordinal
$\alpha$ let $G_{\alpha}(X)$ be the game of length $\alpha$ played
by two players: point picker and open. At each play of the game
point picker picks a real and open responds with an open set
including the real.  Point picker wins a run of the game if at the
end the open sets chosen cover $X$.  The order of $X$ is the least
ordinal for which point picker has a winning strategy.  What orders
are possible?  Daniel and Gruenhage have examples of order
$\omega n$ assuming CH.

\prob (Komjath, see\cite{stepkom})
Suppose every set of reals of size $\omega_1$ has measure zero.
Then does every $\omega_1$ union of lines have planar measure zero?
(Dually) Suppose the real line is union of $\omega_1$ measure zero sets.
Then does there exists $\omega_1$ measure zero subsets of the plane
such that every line is contained in one of them?

\prob (Zhou 3-93)
Does every set of size $\omega_1$ is a Q-set imply that
${\goth  p} > \omega_1$.
For $\gamma$-sets it is true.

\topic{Quasiorder theory}

\prob Is there a Borel version of Fraisse's conjecture?
Are the Borel linear orderings well-quasiordered under embedding?
\answer{Yes, Louveau and St-Raymond \cite{louv} assuming large parts of AD.}

\prob (Laver) Is it consistent that the set of Aronszajn trees is
well-quasi-ordered under embeddability? See Laver \cite{lav}
and Corominas \cite{coromin}.

\prob (Kunen) Is the set of all better-quasi-ordered binary relations
on $\omega$ a proper $\Pi^1_2$ set?
\answer{Yes, Marcone \cite{marco}}

\prob Suppose $(Q,\leq)$ is a recursive quasi-order. Is it true that
Q is BQO iff $Q^{<\omega_1^{ck}}$ is WQO?

\prob (Kunen) Suppose $(Q,\leq)$ is a recursive well quasi-order. Does
$Q^{\omega}/\equiv$ have a recursive presentation?

\prob Suppose every set is Ramsey and $f:\infsets\mapsto\mbox{ORD}$.
Then does there exist $X\in\infsets$ such that the image
of $[X]^\omega$ under $f$ is
countable?  See Louveau-Simpson \cite{simp} and
Aniszcyk-Frankiewicz-Plewik \cite{afp}.

\prob Is finite graphs under homeomorphic embedding WQO?

\prob Is the witness lemma true for LIN(Q) or TREE(Q)?

\prob Is there an $\omega_1$-descending sequence of countable posets
(under embedding) each of which is the union of two chains?
(Kunen, Miller: There is an $\omega_1$-descending sequence of
countable posets.  Kunen: There is an infinite antichain of finite
posets each of which is the union of two chains.)

\prob Is there a parameterized version of Carlson's theorem?
See Carlson \cite{carl} and Pawlikowski \cite{pawpara}.

\topic{not AC}

\prob (M. Bell) Does ZF imply that for every family of nonempty sets
there exists a function assigning to each set in the family
a compact Hausdorff topology? ( Motivation: AC is equivalent
to  this principal
plus every Tychonov product of compact Hausdorff spaces is compact.)

\answer{No, Todorcevic \cite{todlet} Just \cite{just13-1} }

\prob (Dow 88) Does
Stone's theorem on metric spaces
(every metric space is paracompact) require AC?
It is known that ZF implies that $\omega_1$ with the
order topology is not metrizable.
\answer{Yes, Watson see \cite{good}.}

\topicx{Recursion theory}

\prob Does there exist a non-trivial automorphism of the Turing degrees?
(Re degrees?)
\answer{Yes for both? announced by Cooper \cite{cooper}}

\prob (jockusch) Does there exists a DNR of minimal Turing degree?
(DNR means diagonally non recursive: $f\in\bsp$ and
for all $e\in\omega$,
$\;\;f(e)\not=\{e\}(e)$.)

\topic{Miscelleneous}

\prob (Sierpi\'{n}ski) Is there a Borel subset of the plane which
meets every line in exactly two points?
(Mauldin) Must  such a set be zero dimensional?
\answer{Davies has shown such a set cannot be $\Sigma^0_2$ and Mauldin
has shown such a set must be disconnected.  Miller \cite{m24} showed
that if V=L then there does exist a $\Pi^1_1$ subset of the plane
which meets every line in exactly two points.  Kulesza
\cite{kul} showed that any two point set must be zero dimensional.}

\prob (Cichon) Is it consistent to have that the real line is the
disjoint union of $\omega_2$ meager sets such that every meager
set is contained in a countable union of them?
\answer{No, Brendle \cite{brendcic}}

\prob (Juhasz) Does club imply there exist a Souslin line?
\answer{No, Dzamonja-Shelah \cite{dzsh}}

\prob (Ulam \cite{ulam}) Does there exist a set D dense in the plane such that
the distance between any two points of D is rational?

\prob (Miller \cite{m10})
Suppose the continuum is greater than $\omega_2$, then does there
exists a set of reals of cardinality the continuum which cannot be
mapped continuously onto the unit interval?

\prob Is it consistent that there exists $x\in\csp$ such
that $V=L[x]\neq L$ and a continuous onto function
$f:L\cap\csp\mapsto V\cap\csp$?

\prob (Price \cite{price}) Is it consistent there is no Cech function?

\prob (Kunen) Does the consistency of an elementary embedding of $M$ into $V$
imply the consistency of a measurable cardinal?

\prob (Erd\"{o}s) Without CH can you partition the plane into countably many
pieces so that no piece contains an isoceles triangle?
See Kunen \cite{kuntri}.

\prob Is there a Borel version of Hall's marriage theorem?  As for
example, the Borel-Dilworth Theorem \cite{BDT}.

\prob (Davies \cite{davies}) Assuming CH for every
$f:\reals^2\mapsto \reals$ there exists $g_n,h_n:\reals\mapsto \reals$
 such that
$$f(x,y)=\Sigma_{n\in\omega}g_n(x)h_n(y)$$
Does this imply CH?

\prob (Mauldin) CH implies that for every $n\geq 3$ there exists
a 1-1 onto function $f:\reals^{n}\mapsto \reals^{n}$ which maps
each circle onto
a curve which is the union of countably many line segments.  Is CH
necessary?

\prob (Kunen)  Can there be a Souslin tree $T\subseteq 2^{\kappa}$ such
that for all $\alpha<\kappa$ the $T_{\alpha}$ contains all except at
most one of the $\alpha$ branches thru $T_{<\alpha}$.
Here $\kappa$ is
the first Mahlo or weakly Mahlo.

\prob (Baumgartner \cite{baum}) Is it consistent that any two
$\omega_2$ dense sets
of reals are order isomorphic?

\prob (S. Kalikow \cite{kal}) For any set $X$ define for $x,y\in X^{\omega}$,
$\;\;x=^*y$ iff for all but finitely many $n\in\omega$, $\;\;
x(n)=y(n)$.  $X$ has the discrete topology and
$X^{\omega}$ the product topology.
Is it consistent that there exists a map
$f:\omega_2^{\omega}\mapsto \csp$ which is continuous
and for every $x,y\in {\omega_2}^{\omega}$, $\;\; x=^*y$ iff $f(x)=^*f(y)$.
(Kalikow: yes for $\omega_1$ in place of $\omega_2$.)
\answer{Yes, Shelah \cite{shel4}.}

\prob (unknown 1-92) According to Erdos,
Sylvestor proved that given finitely many
points $F$ in the plane not all collinear, there exists a line
$L$ which meets $F$ in exactly two points.   $F={\Bbb Z}\times {\Bbb Z}$
is an obvious infinite counterexample.   Does there exists
a counterexample which is a convergent sequence?  countable compact
set?

\begin{flushright}
                                       Arnold W. Miller\\
                                University of Wisconsin\\
              Department of Mathematics, Van Vleck Hall\\
                                      480 Lincoln Drive\\
                                      Madison, WI 53706\\
                                   miller@math.wisc.edu\\
                              http://www.math.wisc.edu/$\sim$miller\\
\end{flushright}


\begin{thebibliography}{99}

\bibitem{afp}
B.Aniszcyk, R.Frankiewicz, and S.Plewik,
Remarks on (s) and Ramsey-measurable function,
Bulletin de L'Academie Polonaise des sciences, 35(1987),
479-485.

\bibitem{baldowk}
Z.Balogh, A small dowker space in ZFC, preprint 6-95. \where

\bibitem{slalom}
T.Bartoszynski, Combinatorial aspects of measure and category,
Fundamenta Mathematicae, 127(1987), 225-239.

\bibitem{bj}
T.Bartoszynski, H.Judah, On Sierpinski sets, Proceedings of the
American Mathmatical Society, 108(1990), 507-512.

\bibitem{barshel}
T.Bartoszynski, S.Shelah, consistent to have union of strong first
category not strong first category,  announced 95. \where

\bibitem{baum}
J.Baumgartner, All $\aleph_1$-dense sets of reals can be isomorphic,
Fundamenta Mathematicae, 79(1973), 101-106.

\bibitem{beckanal}
H.Becker, Analytic sets from the point of view of compact
sets, Mathematical Proceedings of the Cambridge Philosophical
Society, 99(1986), 1-4.

\bibitem{beck}
H.Becker, The topological Vaught's conjecture and minimal
counterexamples, Journal of Symbolic Logic, 59(1994), 757-784.

\bibitem{becker}
H.Becker, Solution to a particularly interesting problem
of Arnie Miller, handwritten 1-94. \where

\bibitem{brendcic}
J.Brendle, Nicely generated and chaotic ideals, preprint $<95$. \where

\bibitem{bur}
M.Burke, A theorem of Friedman on rectangle inclusion and
its consequences, note dated March 7,1991.

\bibitem{cab}
Cabal Seminar 76-77,77-79,79-81,81-85, Lecture Notes in Mathematics,
689,839,1019,1333, Springer-Verlag.

\bibitem{caic}
X.Caicedo, Independent sets of axioms in $L_{\kappa\alpha}$,
Canadian Mathematical Bulletin, 24(1981), 219-223.

\bibitem{carl}
T.Carlson, Some unifying principles in Ramsey Theory, Disc. Math. 68
(1988) 117-169.

\bibitem{cg}
K.Cielsielski, F.Galvin, Cylinder problem, Fundamenta Mathematicae,
127(1987), 171-176.

\bibitem{cooper}
B.Cooper, Automorphism of the degrees, announced 95. \where

\bibitem{cor}
P.Corazza, Ramsey sets, the Ramsey ideal, and other classes
over $\reals$, Journal of Symbolic Logic, 57(1992), 1441-1468.

\bibitem{coromin}
E.Corominas, On better quasiordering countable trees,
Discrete Mathmatics, 53(1985), 35-53.

\bibitem{davies}
R.O.Davies, Representation of functions of two variables as sums of
rectangular functions I, Fundamenta Mathematicae, 85(1974), 177-183.

\bibitem{vandou}
E.van Douwen, The integers and topology,
in {\bf Handbook of set theoretic topology}, ed by K.Kunen and
J.Vaughan, North-Holland, (1984), 111-167.

\bibitem{dzsh}
M.Dzamonja, S.Shelah, $\clubsuit$ does not imply the
existence of a Suslin tree, preprint 4-96, DjSh 604. \where

\bibitem{erdhech}
P.Erdos, S.Hechler, On maximal almost-disjoint families
over singular cardinals, Colloquia Mathematica Societatis
Janos Bolyai, 10(1973), Infinite and finite sets, Keszthely
Hungary, 597-604.

\bibitem{ekm}
P.Erdos, K.Kunen, and R.Mauldin,
Some additive properties of sets of real numbers,
Fundamenta Mathematicae, 113 (1981),187-199.

\bibitem{ef}
A.Ehrenfeucht, G.Fuhrken, On models iwth undefinable
elements, Math. Scand. 28(1971), 325-328.

\bibitem{falc}
K.Falconer, On a problem of Erdos on sequences and measurable sets,
Proceedings of the American Mathematical Society, 90(1984), 77-78.

\bibitem{fleiss}
W.Fleissner,
An axiom for nonseparable Borel theory,
Transactions of the American Mathematical Society,  251(1979), 309-328.

\bibitem{fleisq}
W.Fleissner, Squares of $Q$ sets,
Fundamenta Mathematicae, 118(1983), 223-231.

\bibitem{fran}
R.Frankiewicz, To distinguish topologically the spaces $m^*$,
Bulletin Academie Polonaise Science, 25(1977), 891-893.


\bibitem{m22}
D.Fremlin and A.Miller,  On some properties of Hurewicz, Menger, and
Rothberger, Fundamenta  Mathematicae, 129(1988), 17-33.

\bibitem{fremreal}
D.Fremlin, Real-valued-measurable cardinals, in {\bf Israel
Mathematical Conference Proceedings}, ed by H.Judah,
6(1993), 151-304.

\bibitem{fried}
H.Friedman, One hundred and two problems in mathematical logic,
Journal of Symbolic Logic, 40(1975), 113-129.

\bibitem{fuh}
G.Fuhrken, A model with exactly one undefinable element,
Colloquium Mathmaticum, 19(1968), 183-185.

\bibitem{good}
C.Good, I.J.Tree, Continuing horrors of topology without
choice, Topology and its applications 63(1995), 79-90.

\bibitem{gros}
M.Groszek, $\omega_1^*$ as an initial segment of c-degrees,
Journal of Symbolic Logic, 59(1994), 956-976.

\bibitem{BDT}
L.Harrington, D.Marker, S.Shelah, Borel orderings, Transactions
of the American Mathematical Society, 310(1988), 293-302.

\bibitem{hech}
S.Hechler, On the existence of certain cofinal subsets of $^\omega\omega$,
in {\bf Axiomatic Set Theory}, proceedings of symposia in pure
mathematics, vol 13 part 2, ed by T.Jech, (1974), 155-174.

\bibitem{hjorem}
G.Hjorth, email message 11-94. \where

\bibitem{hjorthclass}
G.Hjorth, An orbit that is exactly $\Sigma^0_{\lambda +1}$,
handwritten note, 2-95. \where

\bibitem{hjorthmin}
G.Hjorth, On $\aleph_1$ many minimal models, handwritten
note, 12-94. \where

\bibitem{hjorthrig}
G.Hjorth, On $\aleph_1$ many rigid models, handwritten
note, 4-96. \where

\bibitem{hjorth1}
G.Hjorth, Universal co-analytic sets, preprint 12-94. \where

\bibitem{jp}
T.Jech, K.Prikry, Ideals over uncountable sets: applications of
almost disjoint functions and generic ultrapowers, Memoirs of
the American Mathematical Society, 18 no 214, (1979).

\bibitem{js}
H.Judah, S.Shelah, Q-sets do not necessarily have strong
measure zero, Proceedings of the American Mathematical Society,
102(1988), 681-683.

\bibitem{just13-1}
W.Just, email 12-93. \where

\bibitem{m34}
W.Just, A.Miller, M.Scheepers, and P.Szeptycki,
The combinatorics of open covers (II), preprint 1995. \where

\bibitem{kal}
S.Kalikow, Sequences of reals to sequences of zeros and ones,
Proceedings of the American Mathematical Society, 108(1990), 833-837.

\bibitem{kan}
A.Kanamori, {\bf The higher infinite}, Springer-Verlag, (1994).

\bibitem{kecsol}
A.Kechris, S.Solecki, Approximation of analytic by Borel sets
and definable countable chain conditions, preprint 95. \where

\bibitem{kojshe}
M.Kojman, S.Shelah, A ZFC Dowker space in $\aleph_{\omega+1}$:
an application of pcf theory to topology, preprint 1996. \where

\bibitem{kul}
J.Kulesza, A two point set must be zero dimensional,
Proceedings of the American Mathematical Society,
116(1992), 551-553.

\bibitem{m17}
K.Kunen and A.Miller, Borel and projective sets from the point of view of
compact sets, Mathematical Proceedings of the Cambridge Philosophical
Society, 94(1983), 399-409.

\bibitem{kun}
K.Kunen, Random and Cohen reals, in Handbook of
Set Theoretic Topology, North-Holland, (1984), 887-911.

\bibitem{kuntri}
K.Kunen,
Partitioning Euclidean space,
Mathematical Proceedings of the  Cambridge Philosophical Society,
102(1987), 379-383.

\bibitem{kunma}
K.Kunen, Where MA first fails, Journal of Symbolic
Logic, 53(1988), 429-433.

\bibitem{kur}
K.Kuratowski, On a topological problem connected
with the Cantor-Berstein theorem,
Fundamenta Mathematicae, 37(1950), 213-216

\bibitem{lav}
R.Laver, Better-quasi-orderings and a class of trees,
in Studies in Foundations and Combinatorics, Advances in Mathematics
Supplementary Studies  (1978) 31-48.

\bibitem{m25}
S.Leth, J.Keisler, K.Kunen and A.Miller,  Descriptive set theory
on a hyperfinite set, Journal of Symbolic Logic, 54(1989), 1167-1180.

\bibitem{simp}
A.Louveau and S.Simpson, A separable image theorem for ramsey
mappings, Bulletin de L'Academie Polonaise des sciences, 30(1982),
105-108.

\bibitem{louv}
A.Louveau and J. Saint Raymond, On the quasi-ordering of
Borel linear orders under embeddability, Journal of Symbolic
Logic 55(1990), 537-560.

\bibitem{m18}
L.Manevitz and A.Miller,  Lindel\"{o}f models of the reals:
solution to a problem
of Sikorski,  Israel Journal of Mathematics, 45(1983), 209-218.

\bibitem{marco}
A.Marcone, The set of better-quasi-orderings is $\Pi^1_2$,
Mathematical Logic Quarterly, 41(1995), 373-383.

\bibitem{m1}
A.Miller, On the length of Borel hierarchies, Annals of Math Logic,
16(1979), 233-267.

\bibitem{m2}
A.Miller, There are no Q-points in Laver's model for the Borel conjecture,
Proceedings of the American
Mathematical Society, 78(1980), 103-106.

\bibitem{m9}
A.Miller, Generic Souslin sets, Pacific Journal of Mathematics,
97(1981), 171-181.

\bibitem{m10}
A.Miller, Mapping a set of reals onto the reals, Journal of Symbolic Logic,
48(1983), 575-584.

\bibitem{m11}
A.Miller, A characterization of the least cardinal for which the
Baire category theorem fails, Proceedings of the
American Mathematical Society, 86(1982), 498-502.

\bibitem{m13}
A.Miller, The Borel classification of the isomorphism class of a
countable model,
Notre Dame Journal of Formal Logic, 24(1983), 22-34.

\bibitem{m24}
A.Miller, Infinite combinatorics and definability,
Annals of Pure and Applied Mathematical Logic, 41(1989), 179-203.

\bibitem{m27}
A.Miller, Set theoretic properties of Loeb measure,  Journal
of Symbolic Logic, 55(1990), 1022-1036.

\bibitem{m28}
A.Miller, Projective subsets of separable metric spaces,
Annals of Pure and Applied Logic, 50(1990), 53-69.

\bibitem{m31}
A.Miller, Some interesting problems,  Set Theory of the Reals,
ed Haim Judah, Israel Mathematical Conference Proceedings,
vol 6 (1993), 645-654, American Math Society.

\bibitem{m33}
A.Miller, {\bf Descriptive Set Theory and Forcing: how to prove
theorems about Borel sets the hard way},
Lecture Notes in Logic 4(1995), Springer-Verlag.

\bibitem{him}
H.I.Miller, Some results connected with a problem of
Erdos II, Proceedings of the American Mathematical Society,
75(1979), 265-268.

\bibitem{ostem}
Ostaszewski, email 9-92. \where

\bibitem{paw}
J.Pawlikowski, Every Sierpinski set is strongly meager,
preprint $<95$. \where

\bibitem{pawpara}
J.Pawlikowski, Parametrized Ellentuck Theorem,
Topology and its Applications, 37(1990), 65-73.

\bibitem{price}
R.Price, On a problem of \v{C}ech, Topology and Its Applications,
14(1982), 319-329.

\bibitem{m20}
K.Prikry and A.Miller,   When the continuum has cofinality $\omega_1$,
Pacific Journal of Mathematics, 115(1984), 399-407.


\bibitem{shel3}
S.Shelah, {\bf Proper Forcing}, Lecture Notes in Math, 940 (1982),
Springer-Verlag.

\bibitem{shelboulder}
S.Shelah, On cardinal invariants of the continuum,
Collection: Axiomatic set theory Boulder, Colo., 31(1983), 183-208
Contemp. Math., AMS.

\bibitem{shel4}
S.Shelah, On a problem of Steve Kalikow,
preprint (1995) Sh590. \where

\bibitem{shel2}
S.Shelah, Covering of the null ideal may have countable
cofinality, preprint 6-95, Sh592. \where

\bibitem{shel1}
S.Shelah, Embedding Cohen algebras using pcf theory, preprint 7-95,
Sh595. \where

\bibitem{sier}
W.Sierpinski, {\bf Hypothese du Continu}, Monografje Matematyczme,
Warszawa-Lwow, 1934.

\bibitem{sier70}
W.Sierpinski, Problem 70, Fundamenta Mathematicae, 26(1936), 334.

\bibitem{step}
J.Steprans, Almost disjoint families of paths in lattice grids,
Topology Proceedings 16(1991) 185-200.

\bibitem{stepkom}
J.Steprans, Cardinal invariants associated with Hausdorff
capacities, preprint 6-94. \where

\bibitem{solec}
S.Solecki, Covering analytic sets by families of closed sets,
Journal of Symbolic Logic,
59(1994), 1022-1031.

\bibitem{tod}
S.Todorcevic, Remarks on Martin's axiom and the continuum hypothesis,
Canadian Journal of Mathematics, 43 (1991), 832-851.

\bibitem{todlet}
S.Todorcevic, handwritten note  8-91. \where

\bibitem{truss}
J.Truss, Sets having calibre $\aleph_1$, Logic Colloquium 76, North-Holland,
Amsterdam, 1977, 595-612.

\bibitem{ulam}
S.Ulam, {\bf Problems in Modern Mathematics},
Wiley, New York, 1964.

\bibitem{vw}
B.Velikovic, H.Woodin, Complexity of reals of inner models
of set theory, preprint 95. \where


\end{thebibliography}
\end{document}